\documentclass{gtpart}
\usepackage[utf8]{inputenc} 
\usepackage[T1]{fontenc}    
\usepackage{hyperref}       
\usepackage{url}            
\usepackage{booktabs}       
\usepackage{amsfonts}       
\usepackage{nicefrac}       
\usepackage{microtype}      

\usepackage{marvosym}
\usepackage{amsmath}
\usepackage{amsgen}
\usepackage{amsbsy}
\usepackage{amsopn}
\usepackage{amsthm}
\usepackage{amscd}
\usepackage{amsfonts}
\usepackage{amssymb}
\usepackage{mathrsfs}
\usepackage{hyperref}
\usepackage{ifsym}
\usepackage[usenames,dvipsnames]{color}

\usepackage{todonotes}
\usepackage{tikz-cd}
\usepackage{epigraph}

\usepackage[usestackEOL]{stackengine}
\renewcommand{\epigraphsize}{\small}
\newcommand{\mytextformat}{\epigraphsize\itshape}
\newcommand{\mysourceformat}{\epigraphsize\scshape}

\let\originalepigraph\epigraph 
\renewcommand\epigraph[2]{%
  \setbox0=\hbox{\stackon{\textit{\mytextformat\Longstack{#1}}}%
    {\mysourceformat\scshape\Longstack{#2}}}%
  \ifdim\wd0>.8\linewidth\wd0=.8\linewidth\fi%
  \setlength{\epigraphwidth}{\wd0}%
  \originalepigraph{\textit{#1}}{\textsc{#2}}%
}

\usepackage{graphicx}
\usepackage{comment}

\usepackage{latexsym}

\newtheorem{Theorem}{Theorem}
\newtheorem{Lemma}[Theorem]{Lemma}

\newtheorem{Corollary}[Theorem]{Corollary}

\newtheorem{Definition}{Definition}
\newtheorem{Example}[Theorem]{Example}

\def\is{\rm{Iso}}

\newcommand{\cout}[1]{}

\renewcommand{\hat}{\widehat}

\begin{document}
\title[TIE convolutions]{Turing approximations, toric isometric embeddings \\ \& manifold convolutions}

\author[P. Su\'arez-Serrato]{P. Su\'arez-Serrato \\  \today} 

\address{Instituto de Matem\'aticas, Universidad Nacional Aut\'onoma de M\'exico UNAM, Mexico City}

\email{pablo@im.unam.mx}

\begin{abstract}
Convolutions are fundamental elements in deep learning architectures. Here, we present a theoretical framework for combining extrinsic and intrinsic approaches to manifold convolution through isometric embeddings into tori. In this way, we define a convolution operator for a manifold of arbitrary topology and dimension. We also explain geometric and topological conditions that make some local definitions of convolutions which rely on translating filters along geodesic paths on a manifold, computationally intractable. A result of Alan Turing from 1938 underscores the need for such a toric isometric embedding approach to achieve a global definition of convolution on computable, finite metric space approximations to a smooth manifold. 
\end{abstract} 

\maketitle

\epigraph{\it We shall not cease from exploration\\And the end of all our exploring \\Will be to arrive where we started \\ And know the place for the first time.}{T. S. Eliot}
 
 \epigraph{\it Measure or expectation, all must be\\Harvested and yielded}{Seamus Heaney}

\pagebreak

\tableofcontents

 
 \section{Introduction} 

In Convolutional Neural Networks (CNNs) \cite{Sch, CNN-LeNet}, the convolution operations allow for the application of a given filter will to each part of a data file (typically images). Then, as the image moves with translations, activations in each network layer respond with similar translations. This {\it equivariance} property, together with pooling, allows each neuron to express the influence of nearby neurons while training. 

A guiding principle of deep learning is the manifold distribution hypothesis \cite{LBH15}. It posits that high-dimensional data concentrate close to a (nonlinear) lower-dimensional manifold.  The field of manifold learning has been proliferating. Recall that given data (such as cloud points in some ${\bf R}^n$) it is possible to construct a manifold of certain smoothness fitted to them \cite{FIKLN18}. Descriptions of such estimators that approximate these data and are manifolds with bounded reach can be consulted \cite{FMN16}. A review of the development of manifold learning is available \cite{MaFu}. 

The choice of data representation strongly affects the performance of machine learning algorithms. In recent years there has been an increasing interest in extending CNNs to arbitrary, non-euclidean manifolds. A significant challenge has been finding a rigorous definition of convolution on manifolds because addition/subtraction is generally not defined for every manifold. 
 
We propose a new way to define convolutions on manifolds by first isometrically embedding the manifold into a high dimensional torus and then extending a continuous function from the isometric image of the manifold to the target torus. Then, the extended function's convolutions on the torus define the convolution of the original pair of functions. This new definition of convolution is global and works for compact Riemannian manifolds of any dimension.

Informally, we highlight that by first embedding a manifold isometrically into a higher dimensional Euclidean space and fixing a box around it, translations along the axes of this box then permit the definition of a convolution operator. Imagine that the manifold is inside a unit cube, periodically copied along all axes in all directions, and then standard CNNs can be defined on top. Carrying this process out with rigor requires specific control of geometric quantities, which we explain below. 

A toric isometric embedding (TIE) provides a geometric context where discretizations of the ambient space can take into account the intrinsic symmetries of the original manifold completely. The advantage of working with isometric (or even almost isometric) embeddings is that they provide the best of intrinsic and extrinsic worlds. The isometric property preserves the intrinsic geometry, while the global toric coordinates of the embedding allow for convolutions, and in general Fourier analysis, to be carried out. 

Assume the compact smooth connected manifold $M$ is embedded isometrically into the $n$--dimensional torus $T^n$. A function $f$ on $M$ extends to a function $\bar{f}$ on $T^{n}$ (see Lemma \ref{lem-conv-mfds}). Let $k$ be a kernel on $M$, and likewise write $\bar{k}$ for its extension in $T^{n}$. 

Our main contribution is the following Theorem/Definition:

\begin{Theorem} \label{thm-mfds-r-treps}

A global convolution operator can be defined on closed orientable smooth manifolds using toric isometric embeddings (TIE).

\begin{Definition}[TIE convolution on manifolds]\label{def-TIE-conv}
A convolution operator between two functions $f$ and $k$ in $M$, called the TIE convolution and denoted by $f\bowtie k$ can be defined as follows:\[
(f \bowtie k)(x) := \int_{{\bf T}^n}\bar{f}(y)\bar{k}(x-y)dy 
\]
\end{Definition}
\end{Theorem}


Observe that the definition of $f\bowtie k$ is subordinate to an embedding of $M$ into $T^n$. In turn, its computational complexity will depend on the embedding dimension $n$. This approach permits the definition of CNNs on datasets whose elements are smooth manifolds for a fixed embedding method. A discretized version of the TIE-convolution is readily available, as we are now working on a torus (equation (\ref{eqn:disc-3d-conv}) shows one example in 3D).

A notable consequence of Theorem \ref{def-TIE-conv} for the field of {\it geometry processing} is that shapes in 3D space admit a 3D TIE convolution $f\bowtie k$, which can be {\bf globally} defined on meshes, 3D point-clouds, and voxel representations, all with arbitrary topology. These representations can are just the embedding of the data into $3$--dimensional space. 

Moreover, in dimension 3, the explicit computation of the {\em reach} of an embedding of a surface into ${\bf R}^3$ can be achieved using the medial axis. 

As an example, consider the canonical embedding of the sphere $S^2$ in ${\bf R}^3$, realized by unit norm vectors. In this case, the coordinate functions are eigenmaps, and they define an embedding. Similarly, a collection of $n$ eigenmaps can define embeddings of a smooth orientable $d$--manifold into ${\bf R}^n$. These approaches are just illustrative examples because, in practice, we want the embeddings to be isometric.

Performing a convolution using a bump function has the effect of blurring an image or a shape. In this work, we are using such a convolution with a geometrically controlled bump function to increase the dimension where we work, using an isometric embedding, which permits a global convolution to be defined. 

In practice, we need to fix an embedding dimension. There are various available strategies for finding isometric embeddings. The embedding target dimension $n$ affects the performance of the neural networks that use a TIE convolution, as the number of weights in a CNN grows polynomially with $n$. An area of opportunity for improvement here lies in finding the optimal embedding dimension for a given dataset or learning task.

Historically, the problem of finding an isometric embedding of a smooth closed Riemannian manifold was first solved by John Nash \cite{Nash1, Nash2}. Other valuable strategies realize embeddings into $\ell^2$ \cite{BBG94}  and recent advances improve on these ideas using heat kernels \cite{Portegies, LinWu}, and eigenvector fields of the connection Laplacian \cite{ WangZhu}.  

The task of finding isometric embeddings has been implemented using eigenvalues of the Laplace-Beltrami operator \cite{BN02}, using KDE and local PCA \cite{MN17}, or by strengthening Whitney embeddings to produce almost isometric embeddings \cite{FIKLN15} (among plenty of others).  Using a Nash type embedding, the embedding dimension $n$ grows quadratically in $d = \dim M$ \cite{V13}.  


Compared to other, sometimes local, methods of defining convolutions, our approach works efficiently on a manifold $M$ of arbitrary topology. Indeed, translating a filter between points depends on moving between the points and then making sense of how the filter changes. Let $g$ be a smooth Riemannian metric on $M$. When considering geodesics between the points for this task, this strategy requires taking an average over such possible geodesics. In practice, some have proposed using regions where there is a unique geodesic between any two points. This approach is well defined locally, in a chart, but not generally, because a single chart may not cover the entire manifold $M$. Thus, the problem of moving a filter has to consider how the filter changes as it moves along different geodesic paths, and therefore, this includes having to average over the possible geodesics.

The study of the function $C(x,y,\ell)$ that counts the number of geodesics of length at most $\ell$ between $x$ and $y$ has a long history---started by Serre \cite{Serre}---and it is known to have profound connections to the topology and geometry of the underlying manifold \cite{Paternain-book}. 

Recall that an algorithm is {\it efficient} if it runs in polynomial time. 

Topological restrictions to efficient algorithms for computing $C(x,y,\ell)$ on surfaces and $3$--manifolds first come in the form of the growth type of the {\it fundamental group}. In particular, for surfaces and $3$--manifolds, we have the following results:

\begin{Theorem}\label{Thm:sface-intractable-filters}
Let $\Sigma $ be a compact orientable connected surface of genus$>1$. Then for any smooth Riemannian metric $g$ on $\Sigma$ the strategy of averaging filters translated over geodesics between pairs of points $x$ and $y$ on $\Sigma$ is not efficient.  
\end{Theorem}

\begin{Theorem}\label{Thm:3mfds-intractable-filters}
Let $Y$ be a smooth Riemannian $3$-manifold that is neither homeomorphic to a geometric manifold modelled on one of the Thurston geometries ${\bf S}^3, {\bf S}^2\times {\bf R}, {\bf E}^3, {\bf Nil}$, or homeomorphic to a connected sum $L(2,1)\# L(2,1)$ of a lens space $L(2,1)$ whose fundamental group has order $2$, with itself.
Then for any smooth Riemannian metric $g$ on $Y$, the strategy of averaging filters translated over geodesics between pairs of points $x$ and $y$ on $Y$ is not efficient.  
\end{Theorem}

These obstructions highlight the merits of TIE convolutions over other methods. The manifolds left out by Theorem \ref{Thm:3mfds-intractable-filters} are precisely those whose fundamental group has polynomial growth. Thus, in principle, there could be efficient algorithms for computing the geodesic path counting function on these manifolds. 

In terms of homology a well known result of M. Gromov \cite{Gromov} (see also\cite{Paternain-book}) bounds $C(x,y,\ell)$ by the Betti numbers of $M$. Even when the fundamental group is trivial, rational homotopy theory has established when the function $C(x,y,\ell)$ grows exponentially in $\ell$, because it is also bounded below by the growth of the rational homotopy groups of $M$ \cite{FelixHalperinThomas}. For example already in dimension four, the complex plane blown up at three points, ${\bf C}{\rm P} \# 3 \overline{{\bf C}{\rm P}}$, is simply connected and the geodesic counting function $C(x,y,\ell)$ of any smooth Riemannian metric on ${\bf C}{\rm P} \# 3 \overline{{\bf C}{\rm P}}$ has exponential growth. Thus rendering strategies for defining convolutions that rely on translating along geodesics intractable on such a manifold. This phenomenon is explained rigorously by our next result. Recall that a simply connected manifold $M$ is said to be {\it rationally elliptic} if the total rational homotopy $\pi_{\ast}(M)\otimes {\bf Q}$ is finite-dimensional. This means $\pi_{k}(M)\otimes {\bf Q} = 0$ for all $k>k_0$ for some positive integer $k_0$. The manifold $M$ is called {\it rationally hyperbolic} if it is not rationally elliptic.

\begin{Theorem}\label{Thm:rat-hyp-intractable-filters}
Let $M$ be a smooth, closed, simply connected, rationally hyperbolic $n$--manifold, $n\geq 4$. Then for any smooth Riemannian metric $g$ on $M$ the strategy of averaging filters translated over geodesics between pairs of points $x$ and $y$ on $M$ is not efficient.  
\end{Theorem}

This last result leads to the question of how prevalent rationally hyperbolic manifolds are within all manifolds. This kind of intractability is generic for the choice of a Riemannian metric, in the following sense:

\begin{Theorem}\label{Thm:htop-intractable-filters}
Let $M$ be a closed, smooth, $n$--manifold, $n\geq 2$. Then, the strategy of averaging filters translated over geodesics between pairs of points $x$ and $y$ in $M$ is not efficient for a set of $C^\infty$ Riemannian metrics $g$ on $M$ that is open and dense in the space of $C^\infty$ Riemannian metrics equipped with the $C^\infty$ topology. 
\end{Theorem}

The proof uses a profound result by Contreras, which guarantees that Riemannian metrics of positive topological entropy are generic for any smooth manifold \cite{Contreras}. 

These theorems exhibit the potential intractability of some local strategies' topological and dynamical properties to define convolutions. Moreover, they highlight the question of when a manifold admits a metric with $h=0$. In the case of $T^2$, geometric descriptions of such metrics use bands that bound lifts of geodesics to the universal covering \cite{GK}. In general, however, this challenging research direction is very much wide open. 

To end positively, we will now recall a result of Alan Turing, which demonstrates that the only computationally reasonable approach to defining global convolutions uses tori, such as in the TIE-convolution \ref{def-TIE-conv} we defined above.

First, let us go back to the manifold hypothesis and assume that we have data $D$ approximated by a connected smooth manifold $M$. Observe that in order to define a global, computable, convolution operator on $M$, we must also assume the following:

\begin{enumerate}
\item A global group operation can be defined on $M$, in such a way that a convolution may be defined.
\item $M$ can be approximated by a finite metric space $S$ (so that it is computable).
\end{enumerate}

The first property implies that $M$ is a Lie group. As early as 1938 Turing knew that a Lie group that can be approximated by a finite metric space is compact and Abelian \cite[Theorem 2]{Turing38}. Therefore, if we assume both conditions hold, we have shown that our connected $n$--manifold $M$ is a torus $T^{n}$:

\begin{Corollary}[Turing approximations]
A connected $n$--manifold that admits a global convolution operation finitely approximable by a finite metric space is an $n$--torus.
\end{Corollary}

 




Section \ref{sec-reach} reviews the notion of reach and recalls the use of manifolds to approximate data. Section \ref{sec-amenable} explains the notions of group growth and growth of geodesic counting functions, used in the obstructions to efficient strategies mentioned above. Section \ref{sec-toricconv} discusses convolutions on tori. The proofs are found in section \ref{sec-proofs}, and section \ref{sec-conclusions} contains conclusions and suggestions of future work. 


\section{Reach}\label{sec-reach}

Let $M$ be a compact, boundaryless, oriented, smooth and connected manifold. Assume that $M$ is embedded inside the smooth manifold $N$. 

\subsection{Like a rolling stone}

One way to picture the meaning of the reach, before delving into a formal definition, is to imagine a ball $B$ of radius $r$ and dimension equal to dim$\,N$ rolling over $M$, only allowing $B$ to touch $M$ at a single point as it rolls. The largest radius $r_{{\rm max}}$ that satisfies this single contact point condition is the critical radius, or {\it reach}.

Let us now recall the geometric notion of reach, first introduced by Federer \cite{F59} (see also \cite{GuijarroWilhelm}). The reach of a manifold’s embedding measures how it departs from being convex. It takes into account local curvature and global topology. A first definition is: 

\begin{Definition}[reach of a manifold]
The reach $\rho$ of an embedded manifold $M$ is the largest number such that any point at a distance less than $\rho$ from $M$ has a unique nearest point on $M$.
\end{Definition}

Observe that meshes, point clouds, voxel data, or more generally manifolds presented through a locally finite discretization all have a positive reach. Precisely estimating the reach of a data set, or a random manifold, is an active topic of research \cite{AKCRW19, AKTW18}. 

\begin{Example} The reach of a circle of radius $r$ in the Euclidean plane is exactly $r$. Likewise, the reach of the round sphere $S^{n}$ of radius $r$, embedded in ${\bf R}^{n+1}$ in a canonical way, is also $r$. It is well known that smoothly embedded manifolds have positive reach \cite{Thale}.
\end{Example}

\begin{Example}\label{Ex:MNIST} Consider MNIST. Numbers are idealized as (essentially) $1$D objects. However, the data in MNIST already have a neighborhood about each digit. A number $1$ is represented as embedded in a square $Q$ in ${\bf R}^2$, together with a small neighborhood $\nu(e(1))$ of the image $e(1)\subset Q$. The thickness of the trace, seen as $\nu(e(1))$, is smaller than the reach of $e(1)$. Thus, bounding the thickness of the trace is an important feature for numbers like $0$, $6$, or $9$. If the thickness of the trace were not smaller than the reach, we would see a disk and not a circle. This example illustrates how the reach is an important geometric feature of the representation of digits provided by the MNIST dataset.
\end{Example}

\begin{figure}
\centering
\includegraphics[width = .7 \textwidth]{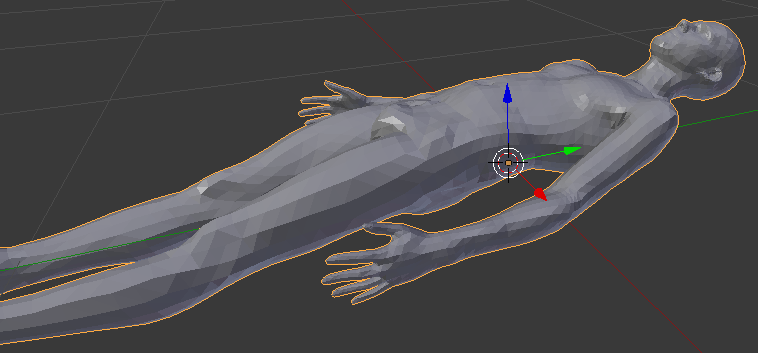}
\caption{The first file from the FAUST dataset, naturally embedded in ${\bf R}^3$. Image courtesy of Eduardo Velazquez Richards. }
\label{im-FAUST1}
\end{figure}

\begin{figure}
\centering
\includegraphics[width = .7 \textwidth]{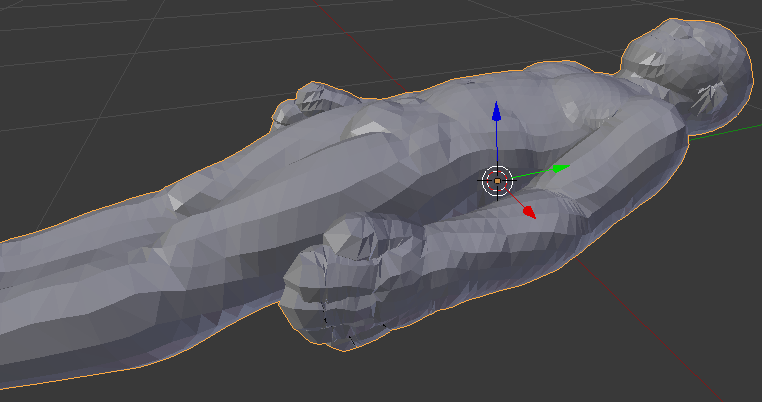}
\caption{The image shows a tubular neighborhood of the first file from the FAUST dataset, which has a width that goes beyond the reach of the embedding. Notice how the fingers have merged. Image courtesy of Eduardo Velazquez Richards.}
\label{im-FAUST1-overreach}
\end{figure}

\begin{figure}
\centering
\includegraphics[width = .7 \textwidth]{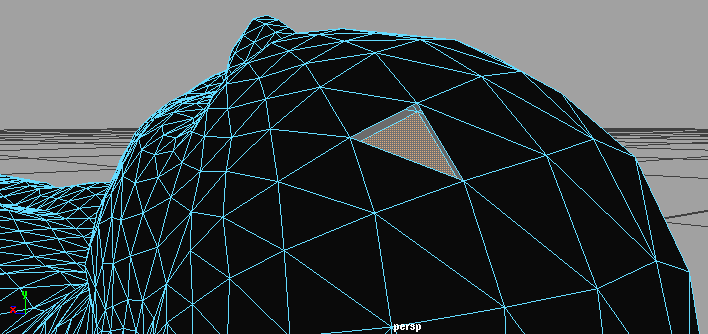}
\caption{Here we see a closer look into a tubular neighborhood of the first file from the FAUST dataset, which has a width below the reach of the embedding. Notice that some cells were removed to emphasize this point. Image courtesy of Eduardo Velazquez Richards.}
\label{im-FAUST1-reachdetail}
\end{figure}

\begin{Example}\label{Ex:FAUST} Consider the FAUST dataset of human 3D poses \cite{FAUST}. In Figure \ref{im-FAUST1} we observe the first mesh in the FAUST dataset. It is already embedded in ${\bf R}^3$. In Figure \ref{im-FAUST1-overreach} a tubular neighborhood of the first mesh in FAUST of thickness beyond the reach is shown. Figure \ref{im-FAUST1-reachdetail} zooms in to show a neighborhood of the first FAUST mesh whose thickness is below the reach of the embedding. 
\end{Example}

\subsection{Rigorous reach}

Consider a point $x$ in $M$ and a unit element $v$ in $T_{x}N$, of the tangent space to $N$ at $x$. Let $\gamma_{x,v}$ be the unit speed geodesic in $N$ based at $x$ in the direction of $v$. The reach is intimately related to the largest balls radii around the origins in $T_{x}N$, for $x\in M$, such that all their exponential maps are diffeomorphisms. 

For $x$ in $M$ and $A\subset M$ define the distance from $x$ to $A$ as:
\[
d_{M}(x,A) \, := \, \inf\limits_{y\in A} d_{M}(x,y)
\]
Here $d_{M}(x,y)$ is the geodesic distance between $x$ and $y$, inside $M$. Likewise, define $d_{N}(x,A)$, for $x$ in $N$, $A\subset N$ and $d_{N}$ the geodesic distance inside $N$.

Then, the {\it local reach} of $M$ in $N$, in a unit direction $v\in T_{x}N$ is defined by \cite{AKTW18}:
\[
\rho_{{\rm loc}}(x,v)\, :=\, \sup \{\, p\, : \, d_{N}(\exp_{x}^{N}(p, v),\, M)=p\}
\]
 So that if $p$ is larger than $\rho_{{\rm loc}}(x,v)$, then there exists a point $y\neq x$ in $M$ which is closer to $\exp_{x}^{N}(p, v)$ than $x$ is. 
 
 We define the local reach of $M\subset N$ at $x$ as:
 \[
 \rho(x)\, := \, \inf\limits_{v\in U\nu_{x}M} \rho_{{\rm loc}}(x,v)
 \]
 Here, $U\nu_{x}M$ means the intersection of the normal space at $x$ with the unit tangent space of $N$ at $x$. 
Finally, taking an infimum, we establish the (global) reach $\rho(M)$ of $M$ to be:
 \[
 \rho(M)\, := \, \inf\limits_{x\in M} \rho(x)
 \]
 
\subsection{Using the reach}

The reach of an embedding can be used to extend functions inside the ambient manifold, as the next result illustrates. 

\begin{Lemma}\label{lem-conv-mfds}
Let $f:M\to {\bf R}$ be a continuous function. Assume $M$ is isometrically embedded in $T^n$, of positive reach $\rho$. Then $f$ can be extended to $\bar{f}:T^{n}\to {\bf R}$.
\end{Lemma}
Tietze's extension theorem implies this result. However, most proofs show only existence. We will now give a constructive argument, which is possible to implement in a numerical scheme.

\begin{proof} Let $\nu(M)$ be a tubular neighborhood of $M$ inside $T^{n}$ of radius equal to $\rho/2$, half the reach. Observe that inside $\nu(M)$ there exist parallel copies (that do not intersect) of $M$, filling out $\nu(M)$ along the normal direction. Therefore, $f$ extends to a function $\breve{f}: \nu(M) \to {\bf R}$. Extend $f$ radially from $M$ by using the same value as the central copy of $M$ on the parallel copies inside $\nu(M)$. 

Notice that, by the definition of reach, there exists a unique point $y$ in $\nu(M)$ at distance $\delta<\rho$ to $x$ in $M$. Thus, the value $f(x)$ can be assigned to this $y$, effectively defining $f$ as the same value on the level sets of the distance function as on $M$. 

Now, multiply this value by a smooth bump function $\beta$, such that $\beta = 1$ on $M$ and $\beta = 0$ on the complement $\nu(M)^{c}$ of $\nu(M)$. Therefore, the following formula is obtained for $y$ in $\nu(M)$:
\[
\breve{f}(y):= \beta(y)\cdot f(x)
\]

Recall that $x$ was the closest point to $y$ in $M$ (unique by definition of the reach). 
Now, extend $\breve{f}$ to a continuous function $\bar{f}: T^{n}\to {\bf R}$, setting $\bar{f}(y)=\breve{f}(y)$, for $y$ in $\nu(M)$, and zero otherwise.
\end{proof}

Notice that positive reach manifolds are homotopy equivalent to metric Vietoris-Rips complexes (of scale parameter below the reach) built on them (cf. Theorem 4.6. in \cite{AM19}). This construction, therefore, completely preserves the topology of the original object. 

\section{Growth}\label{sec-amenable}


\subsection{Group Growth}

\subsubsection{Growth functions}

The growth of a space, such as a manifold or a group, can be registered using functions with the following properties. A function $a: [0, \infty ) \to {\bf R}$ is called a growth function if $a(0) \geq 1$, it is monotonically
increasing, and $a$ is submultiplicative, that is, for all $r, s \geq 0$ and a constant $C_{a}$,
\[
 a(r + s) \leq C_{a} a(r)a(s).
 \]

For $\alpha > 0 $, the growth function $(r+1)^{\alpha}$ is said to be of {\em polynomial} growth, and the growth function $e^{\alpha r}$ is said to be of {\em exponential} growth.

\subsubsection{Growth of generating sets}
Let $S$ be a finite and symmetric generating set of a group $\Gamma$. Denote by 
$N_{S}(m)$ be the number of elements of $\Gamma$ that can be expressed as a word of
length at most $m \in {\bf N }\setminus \{ 0 \}$ in S. Then $N_{S}$ is monotonically increasing, $N_{S}(0) = 1$ and
$N_{S}(m + n) \leq N_{S}(m)N_{S}(n)$.

 Setting $C_{a}=N_{S}(1)$ makes $N_{S}(\lfloor r \rfloor)$ a growth function, because $\lfloor r + s \rfloor \geq \lfloor r \rfloor + \lfloor s \rfloor + 1$.

 \subsubsection{Growth types of groups}
 
 The growth {\em type} of a group is independent of the choice of generating set $S$. So, for example, if the growth of a particular generating set is exponential, then the growth of any other generating set will also be exponential. Therefore the concepts of growths of exponential or polynomial type are well defined for groups. 
 
 Consider the following examples that are relevant for computer vision and geometry processing.
 
 \begin{Example}[The $2$--sphere $S^2$] 
 As $\pi_1(S^2)$ is the trivial group, there is no growth because the identity is the only element. However, this example is not irrelevant. Some geometric deep learning models have been successfully implemented on data that lies on the surface of $S^2$ because its fundamental group does not impose computational costs. 
 \end{Example}
 
 \begin{Example}[The $2$--torus $T^2$]
 The fundamental group of the $2$--torus is isomorphic to ${\bf Z}\oplus{\bf Z}$. Therefore the growth of $\pi_1(T^2)$ is bounded above by a polynomial of order two. Informally, as $T^2$ lifts to its universal covering ${\bf R}^2$, its fundamental group is isomorphic to the integer lattice ${\bf Z}\oplus{\bf Z}$. A ball centered at the origin will contain several integer lattice elements. These grow quadratically as the radius expands. 
 \end{Example}

\subsection{Growth of geodesic arcs}

We will now review some notions connected to the function $C(x,y,\ell)$ and which we will use in the proofs of our intractability results mentioned in the introduction. We refer interested readers to G.P. Paternain's book on geodesic flows for a comprehensive treatment of these concepts \cite{Paternain-book}.

The first one is the {\it topological entropy of the geodesic flow}, which is best understood in a Riemannian context by Ma\~n\'e's formula:

\begin{Theorem}(Ma\~n\'e \cite{Paternain-book})\label{thm:mane}
If a Riemannian metric $g$ on a compact manifold $M$ is of class $C^{\infty}$ then the topological entropy $h$ of the geodesic flow of $g$ equals 
\[ \lim_{\ell \to \infty} \frac{1}{\ell}\int_{M} \log(C(x,y,\ell)) dx\, dy .
\]
\end{Theorem}

These relationships link the dynamics of the geodesic flow and the average rate of growth of geodesics between two points. Thus, if the topological entropy of the geodesic flow is positive, then on average, there are exponentially many geodesic arcs between any two points in the manifold. 

Fix a compact domain $N$ for the action of $\pi_1(M)$ on the universal cover $\widetilde{M}$. Write $a$ for the diameter of $N$. Consider the set $F\subset \pi_1(M)$ of elements $\alpha$ such that $\alpha N \cap N \neq \emptyset$. Let $\nu>0$ denote the exponential growth rate of $\pi_1(M)$ with respect to $F$. 

\begin{Theorem}(Dinaburg \cite{Paternain-book}) \label{thm:dinaburg}
\[
h>\nu /2a.
\]
\end{Theorem}

We are now prepared to state the most relevant implication of exponential group growth for our purposes. It involves both the entropies we just mentioned and is given by the following:

\begin{Lemma}\label{lemma:pi1growth}
Assume the closed smooth Riemannian manifold $M$ has a fundamental group $\pi_1(M)$ whose growth is of exponential type. Then, the function $C(x,y,\ell)$ grows faster than any polynomial in $\ell$.
\end{Lemma}

\begin{proof}
We will show the result holds by contradiction. Assume then that $\pi_1(M)$ has exponential growth type and that $C(x,y,\ell)$ grows at most polynomially in $\ell$. Taking the logarithm, averaging over $M$ by integrating and then taking the limit while dividing by the length $\ell$ we find:
\begin{equation}\label{formula}
\lim_{\ell \to \infty} \frac{1}{\ell}\int_{M\times M} \log(C(x,y,\ell)) dx\, dy\, := h
\end{equation}
Therefore $h=0$ because $C(x,y,\ell)$ grows at most polynomially in $\ell$. By Ma\~n\'e's formula in Theorem \ref{thm:mane}, $h$ equals the topological entropy of the geodesic flow. Moreover, by Dinaburg's Theorem \ref{thm:dinaburg} we know $h>v/2a >0$. This contradiction yields the result. 
\end{proof}

\begin{figure}[htbp]
\centering
\includegraphics[width =  \textwidth]{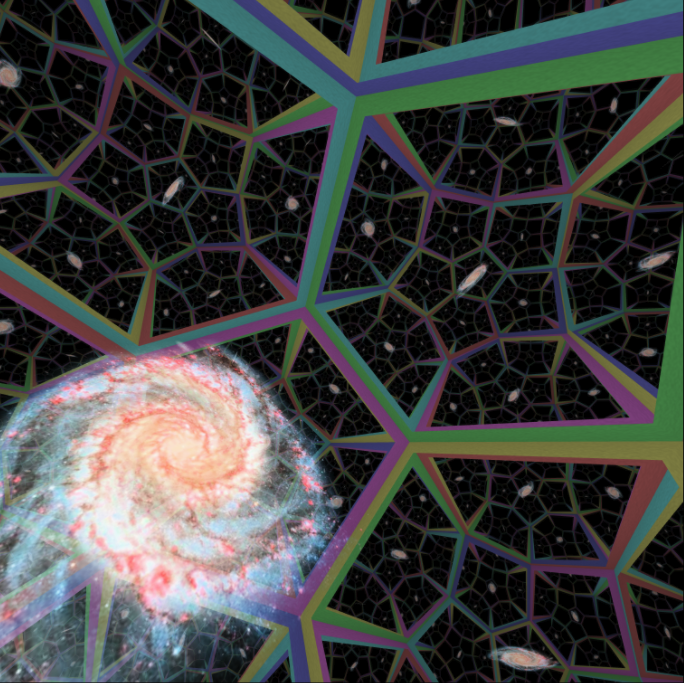}
\caption{A peek inside a negatively curved universe. The number of galaxies we see grows exponentially as we cross more and more fundamental domains. Image generated with the {\it Curved Spaces} software.}
\label{im-hyperbolic-universe}
\end{figure}

\begin{figure}[htbp]
\centering
\includegraphics[width = \textwidth]{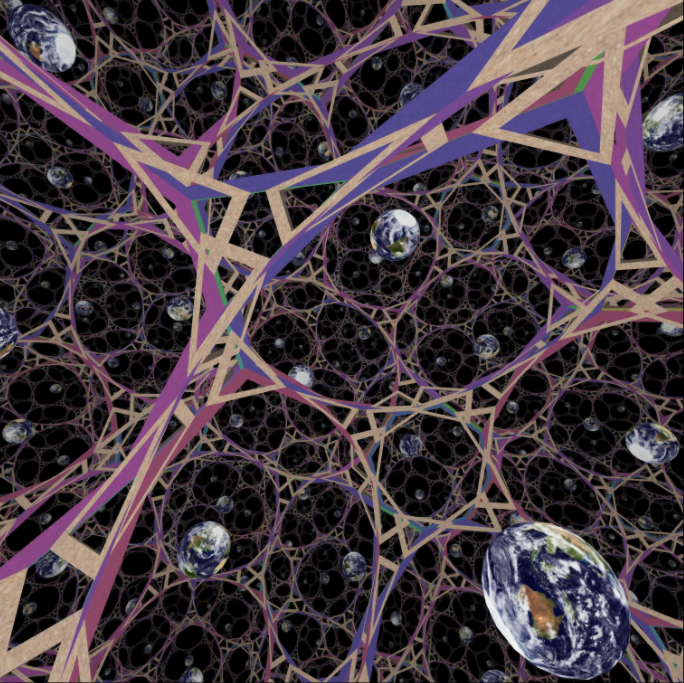}
\caption{A planet living in a hyperbolic $3$-manifold, seen from inside. The copies of the planet that the inhabitants could see proliferate exponentially.  Image generated with the {\it Curved Spaces} software.}
\label{im-hyperbolic-world}
\end{figure}

\begin{Example}\label{ex:exp-growth-neg-curv}
Milnor showed that the fundamental group of a manifold that admits a metric of negative sectional curvature has exponential growth \cite{Milnor}. Thus, for example, surfaces of genus$>2$ and quotients of real hyperbolic spaces. See Figures \ref{im-hyperbolic-universe} and \ref{im-hyperbolic-world}. 
\end{Example}

\subsection{Growth from rational homotopy}
In this subsection, we will assume the manifold $X$ is simply connected.
Maps of higher dimensional spheres into a manifold define the {\it higher homotopy groups}.

\begin{Definition}
The group of homotopy equivalences of maps from $S^{k}$ into the manifold $M$ is denoted by $\pi_{k}(X)$, called the $k-th$ homotopy group of $X$.
\end{Definition}

Tensoring all of these groups with the rational numbers, we find the {\it total rational homotopy} $\pi_{\ast}(X)\otimes {\bf Q}$. This algebraization of the homotopy groups has proven very powerful ever since it was first introduced by Quillen \cite{Q69}, and Sullivan \cite{S77}. For our purposes, it is sufficient to understand when the total rational homotopy of a space is finite-dimensional, as we will see below. We first recall the following:

\begin{Definition}
A manifold is said to be {\it rationally elliptic} if the total rational homotopy $\pi_{\ast}(X)\otimes {\bf Q}$ is finite-dimensional. That is, there exists a positive integer $i_0$ such that for all $i>i_0$, $\pi_{\ast}(X)\otimes {\bf Q}=0$. The manifold $X$ is said to be rationally hyperbolic if it is not rationally elliptic.
\end{Definition}

Let us now review some examples of rationally elliptic and rationally hyperbolic manifolds (see \cite{FelixHalperinThomas, Paternain-book}).

\begin{Example}
A simply connected homogeneous space is rationally elliptic.
\end{Example}

\begin{Example}[Lemma 5.2 \cite{Paternain-book}]
A large enough connected sum of copies of a manifold that is not a rational homology sphere is rationally hyperbolic.
\end{Example}

\begin{Example}[Lemma 5.3 \cite{Paternain-book}]
Let $X$ be any simply connected compact manifold of dimension four or five. If the homology of $X$ satisfies ${\rm dim}\, H_{2}(X;{\bf Q})>2$ then $X$ is rationally hyperbolic.
\end{Example}

The principal consequence of rational hyperbolicity in the context that interests us here is aided by the following result of G.P. Paternain:

\begin{Theorem}[Paternain, Corollary 5.21 \cite{Paternain-book}]\label{thm:rat-hyp-h}
Let $X$ be a rationally hyperbolic manifold. Then for any $C^{\infty}$ Riemannian metric on $X$ the topological entropy of the geodesic flow $h$ satisfies $h>0$.
\end{Theorem}

The arguments that show this last result do not hold for manifolds with an infinite fundamental group. It is thus not clear in that case how to relate the growth of the average number of geodesics between two points to the homology of the loop space and hence to rational hyperbolicity.

\subsection{Growth in generic Riemannian metrics}\label{subsec:generic}

Gonzalo Contreras proved a perturbation lemma for the derivative of geodesic flows in high dimensions. Therefore, a generic metric has a non-trivial basic set in its geodesic flow and has positive topological entropy. 

Recall that, in a topological space, a residual set is the complement of a countable union of nowhere dense sets. In particular, a residual set is dense in its ambient space. This notion allows us to formalize the idea of generic metrics in a set-theoretic way.

One of the relevant results of Contreras is the following:

\begin{Theorem}[Theorem C, \cite{Contreras}]\label{thm::Contreras-ThmC}
Let $G_0$ be the set of $C^4$ Riemannian metrics on $M$ such that, \\(i) The K-jet of the Poincar\'e map of every closed geodesic of $g$ belongs to $Q$; \\(ii) all heteroclinic points of hyperbolic closed geodesics of $g$ are transversal. 

Then;
\begin{enumerate}
\item $G_0$ contains a residual set in $R^{k}(M)$ for all $k\geq 4$.
\item If the geodesic flow of a metric $g\in G_0$ contains hyperbolic periodic orbit then it contains a non-trivial hyperbolic set. In particular $h_{top}(g)>0$.
\end{enumerate}
\end{Theorem}

Recall that, in a topological space, a residual set is the complement of a countable union of nowhere dense sets. A non-empty complete metric space with a non-empty interior is not the countable union of nowhere dense sets.

Let $\mathcal{R}^{k}$ be the space of $C^{k}$ Riemannian metrics on $M$ equipped with the $C^{k}$ topology. Notice that because $\mathcal{R}^{\infty}$ is a complete metric space, residual sets are dense within it.

The next result by Contreras is the most pertinent for our work here: 

\begin{Theorem}[Theorem A, \cite{Contreras}]\label{thm::Contreras-ThmA}
On $M^{n}$, $n\geq 2$, the set of $C^{\infty}$ metrics whose geodesic flow on the unit tangent bundle admits a nontrivial hyperbolic basic set is open and dense in the $C^{2}$ topology.
\end{Theorem}

Therefore, on $M$ the set of $C^{\infty}$ metrics $g$ with $h_{top}(g)>0$ contains an open and dense set in the $C^{2}$ topology. This is the precise sense in which smooth, $C^{\infty}$, {\it generic} metrics have $h_{top}(g)>0$.

\section{Toric Convolutions}\label{sec-toricconv}

\subsection{Properties of CNNs}

In a feed forward neural network, {\it neurons} are arranged in $L+1$ distinct layers, input at $l=0$, and output at $l=L$. Denote by $\chi_{l}$ the index set of the layer and let $V_{l}$ be a vector space. Then we can regard the {\it activation functions} of the network in a given layer as functions $f^{l}:\chi_{l} \to V_{l}$. As:
\[
f^{l}(x) = \sigma(b^{l}(x)+ \sum_{y} w^{l}(x,y)\cdot f^{(l-1)}(y))  
\]
Here each layer can be seen as as affine transformation $\phi_{l}: L_{V_{l-1}}(\chi_{l-1}) \to L_{V_{l}}(\chi_{l})$, followed by $\sigma$.

\begin{Definition}
Let $\chi_{0}, \ldots, \chi_{L}$ be a sequence of index sets $V_{0}, \ldots , V_{L}$ vector spaces, $\phi_{0}, \ldots , \phi_{L}$ affine maps, $\phi_{l}:L_{V_{l-1}}(\chi_{l-1})\to L_{V_{l}}(\chi_{l})$, and $\sigma_{l}:V_{l}\to V_{l}$ activation functions. These data define a multi-layer feed-forward neural-network through the sequence of compositions:
\[
f_0\mapsto f_1 \mapsto f_2 \mapsto \ldots \mapsto f_{l}\,;\; f_{l}(x)= \sigma_{l}(\phi_{l}(f_{l-1}(x)))
\]
\end{Definition}

A CNN is a multi-layer feed-forward neural network with these characteristics.

For example neural nets used in image recognition typically use $\chi_0= [m]\times [m]$, and similarly define $\chi_{i}$. The functions $\phi_{l}$ used are such that:
\[
\phi_{l}(f_{l})
(x_1, x_2)=\sum\limits_{y_1=1}^{m}\sum\limits_{y_2=1}^{m}f_{l-1}(x_1 - y_1 , x_2 - y_2) \chi_{l}(y_1, y_2)
\] 
This function is known as the discrete convolution of $f_{l-1}$ with the filter $\chi$. In practice, the "width" of the filters is small, between 3 and 10 pixels, while the number of layers can be quite deep.

The same idea can be used in higher dimensions. One of the purposes of this note is to contribute to the theoretical underpinnings of specific recent approaches to voxel neural networks that use precisely the approach described next.

For a voxel CNN for objects of $m^3$ voxels, a similar network as for images can be defined. Set $\chi_0= [m]\times [m] \times [m]$, similarly define $\chi_{i}$. Set:
\begin{equation}\label{eqn:disc-3d-conv}
\phi_{l}(f_{l})
(x_1, x_2, x_3)=\sum\limits_{y_1=1}^{m}\sum\limits_{y_2=1}^{m}\sum\limits_{y_3=1}^{m}f_{l-1}(x_1 - y_1 , x_2 - y_2, x_3 - y_3) \chi_{l}(y_1, y_2, y_3)
\end{equation} 

Implementations of this type of architecture already exist \cite{MS15, MNA16, GGGOCA16,ZT18}.  

The main contribution of this paper is to be able to define the {\it minus } sign in the previous two formul\ae\, for general manifolds. This aim is achieved through isometric embeddings, metric thickenings, and extensions of functions. Finally, these techniques define convolution on any smooth, compact, orientable manifold, as demonstrated below.

\subsection{Fourier transforms on tori}

Tori, such as ${{ T}^d :=({\bf R}/{\bf Z})^d}$ , with the standard topology, are examples of compact abelian groups. 

In such abelian groups $G$ there exists a {\it translation}:  ${x \mapsto \tau_x}$ of ${G}$. It acts on continuous functions of compact support $ f \in C_c(G)$. For every ${x \in G}$, the translation operation is given by 
\[
{\tau_x: C_c(G) \rightarrow C_c(G)}\, ;\, \tau_x f(y) := f(y-x). 
\]

On groups, a useful measure is the Haar measure. It is a translation-invariant Radon measure ${\mu}$ on ${G}$, that is ${\mu(E+x) = \mu(E)}$ for all Borel sets ${E \subset G}$, ${x \in G}$. Here ${E+x := \{ y+x: y \in E \}}$ will be called the translation of ${E}$ by ${x}$. 

Therefore integration ${f \mapsto \int_G f\ d\mu}$ with respect to a Haar measure is translation-invariant:

\[
  \int_G f(y-x)\ d\mu(y) = \int_G f(y)\ d\mu(y) 
  \]
or equivalently,
\[
 \displaystyle  \int_G \tau_x f\ d\mu = \int_G f\ d\mu 
 \]
for all ${f \in C_c(G)}$, ${x \in G}$.

On  ${({\bf R}/{\bf Z})^d}$, a Haar (probability) measure is obtained by identifying this torus with ${[0,1)^d}$ as usual and taking the Lebesgue measure. 

A  continuous homomorphism ${\chi: G \rightarrow S^1}$ to the unit circle in ${\bf C}$, ${\chi(x+y)=\chi(x) \chi(y)}$ for all ${x,y \in G}$,  is called a {\it multiplicative character}. A continuous homomorphism ${\xi: G \rightarrow {\bf R}/{\bf Z}} $ is called an {\it additive character} or frequency 
\[
{\xi: x \mapsto \xi \cdot x}, \, {\rm so}\,\, {\xi \cdot (x+y) = \xi \cdot x + \xi \cdot y},
\] 
 for all ${x,y \in G}$.
The Pontryagin dual of ${G}$ is defined to be the set of all frequencies ${\xi}$ and is denoted ${\hat G}$; it is an abelian group. A multiplicative character is called non-trivial if it is not the constant function ${1}$; an additive character is called non-trivial if it is not the constant function ${0}$. For an additive character  ${\xi \in \hat G}$ the function ${x \mapsto e^{2\pi i \xi \cdot x}}$ is a multiplicative character. Conversely, every multiplicative character can be described uniquely from an additive character in this way.

The Pontryagin dual ${\widehat{({\bf R}/{\bf Z})^d}}$ of ${({\bf R}/{\bf Z})^d}$ is isomorphic to  ${{\bf Z}^d}$, this is seen identifying each ${\xi \in {\bf Z}^d}$ with the frequency ${x \mapsto \xi \cdot x}$ given by the dot product.

The formula
\[
 \hat f(\xi) := \int_G f(x) e^{-2\pi i \xi \cdot x}\ d\mu(x),
\]
defines the Fourier transform of an absolutely integrable function ${f \in L^1(G)}$.

It is a linear transformation, with the following bound:
\[
 \sup_{\xi \in \hat G} |\hat f(\xi)| \leq \|f\|_{L^1(G)}
\]

On the one hand, translations are converted into frequency modulations, as,
\[
 \widehat{\tau_{x_0} f}(\xi) = e^{-2\pi i \xi \cdot x_0} \hat f(\xi) 
\]
for any ${f \in L^1(G)}, {x_0 \in G}$, and ${\xi \in \hat G}$. On the other hand, frequency modulations are converted to translations, 
\[
 \widehat{\chi_{\xi_0} f}(\xi) = \hat f(\xi-\xi_0) 
\]
for any ${f \in L^1(G)}$ and ${\xi_0,\xi \in \hat G}$, where ${\chi_{\xi_0}}$ is the multiplicative character 
\[
{\chi_{\xi_0}: x \mapsto e^{2\pi i \xi_0 \cdot x}}.
\]

For ${f,g \in L^1(G)}$, the convolution ${f*g: G \rightarrow {\bf C}}$ is defined by the following equation:
\begin{equation}\label{eqn:toric-conv}
 f*g(x) := \int_G f(y) g(x-y)\ d\mu(y)
\end{equation}
Young’s inequality implies that ${f*g}$ is defined a.e. and is in ${L^1(G)}$. In fact,
\[
 \|f*g\|_{L^1(G)} \leq \|f\|_{L^1(G)} \|g\|_{L^1(G)}.
\]

The convolution operation ${f, g \mapsto f*g}$ is bilinear, continuous, commutative, and associative on ${L^1(G)}$. Therefore, the Banach space ${L^1(G)}$ becomes a commutative Banach algebra with this convolution operation as “multiplication”. Defining  ${f^*(x) := \overline{f(-x)}}$ for all ${f \in L^1(G)}$, this turns ${L^1(G)}$ into a Banach {*}-algebra.

For ${f, g \in L^1(G)}$,
\[
\widehat{f*g}(\xi) = \hat f(\xi) \hat g(\xi) 
\]
for all ${\xi \in \hat G}$, so the Fourier transform converts convolution to a pointwise product.

For  ${G = { T}^d = ({\bf R}/{\bf Z})^d}$, the multiplicative characters ${x \mapsto e^{2\pi i \xi \cdot x}}$ separate points, in the sense that given any two ${x,y \in G}$, there exists a character that takes different values at ${x}$ and at ${y}$. The space of finite linear combinations of multiplicative characters (i.e., the space of trigonometric polynomials) is then an algebra closed under conjugation that separates points and contains the unit ${1}$. Hence, by the Stone-Weierstrass theorem, it is dense in ${C(G)}$ in ${L^2}$.

\section{Proofs}\label{sec-proofs}

\subsection{Proof of Theorem \ref{thm-mfds-r-treps} and Definition \ref{def-TIE-conv}}
\begin{proof}
Let $(M, g)$ be a compact smooth Riemannian manifold of dimension $d\geq 0$. Denote by $e:M\to {\bf R}^n$ an isometric embedding, produced, for example, through a Nash embedding. Thus $e$ is a smooth injective immersion of the form,
\[
\langle \partial_{y_{i}} e ,  \partial_{y_{i}} e \rangle_{{\bf R}^n} = g( e_i, e_j),
\] 
in local coordinates $y_1, \ldots, y_n$ with respect to the Euclidean inner product $\langle \, ,  \, \rangle_{{\bf R}^n}$ in ${\bf R}^n$.

Let $\rho>0$ be the reach of $e(M)$. Define $\nu(M)$ to be a tubular neighborhood of width $\rho /2$ about $e(M)$ inside ${\bf R}^n$. Set ${\rm diam}(e(M))$ as the diameter of the image of $M$ under embedding $e$. 

We'll now put the image of $M$ under the embedding into a box $B$, as follows. Consider an $n$-dimensional parellelepiped $B$ of diameter $2\cdot {\rm diam}(e(M)$ that contains $e(M)$ whose sides are all of the same length $\ell$ and each parallel to a ${\bf R}^n$ coordinate axis $x_i$.

Let $q$ denote the translation by $\ell$ on each coordinate. Identify opposite faces of $B$ by taking the quotient under the action of $q$, to form an $n$-torus ${\bf T}:= B / q $.

In this way we have constructed a map $\tau : M \to T^{n}$, defined by $\tau(x)= q\circ e(x)$, for $x$ in $M$, and obtained an isometric embedding of $M$ into $T^{n}$. 

Recall that the standard convolution operator for real functions $f$ and $k$ on Euclidean ${\bf R}^n$ is defined by:
\[
(f \ast k)(x) := \int_{{\bf R}^n}k(x-y)f(y)dy 
\]

By  Lemma \ref{lem-conv-mfds} the functions $f,k$ now defined on a manifold $M$ extends to a functions $\bar{f},\bar{k}$ on $T^{n}$.

Therefore we may now define, as in equation (\ref{eqn:toric-conv}),
\[
(f \bowtie k)(x) := \int_{{\bf T}^n}\bar{k}(x-y)\bar{f}(y)dy. 
\]
\end{proof}


\subsection{Proof of Theorem \ref{Thm:sface-intractable-filters}}
\begin{proof} We will rely on Lemma \ref{lemma:pi1growth}. Hence we will first show that the fundamental group of $\Sigma $ has exponential type growth. We include the following details to guide interested readers. As we have mentioned already by Milnor's result (see Example \ref{ex:exp-growth-neg-curv}), these groups are known to have exponential growth. To guide the argument, we will divide the proof into two cases. 
{\bf Compact orientable surfaces of genus$=2$:} 
 First consider a compact orientable surface $\Sigma$ of genus $2$. The universal covering of $\Sigma$ to a unit $2$--disk $D$. Assume that $\Sigma$ is given a Riemannian metric of Gaussian curvature equal to $-1$ at every point. This metric naturally lifts to a metric $h$ on $D$. The fundamental group $\pi_1(\Sigma)$ is isomorphic to the group of deck transformations of the covering $D\to \Sigma$, acting by isometries of $h$. Fix a point $o\in D$ and write $\mathcal{O}$ for the orbit of $o$ in $D$ under the action of $\pi_1(\Sigma)$. The set $\mathcal{O}$ is a geometric representative that will help us to calculate the growth of $\pi_1(\Sigma)$. Consider a fundamental octagon $P\subset D$, that under the action of $\pi_1(\Sigma)$ develops into a hyperbolic tessellation $\tau$ of $D$. A dual tesselation $\tau '$ can be found by connecting pairs of points in $\mathcal{O}$ whenever they lie in copies of $P$ that share a side. In this way, the growth of the $1$--skeleton of $\tau '$ gives us the growth of $\pi_1(\Sigma)$. Observe that the hyperbolic geometry of $h$ makes the balls of increasing radius $r$ centered at $o$ contain a number of elements of $\pi_1(\Sigma)$ that grows exponentially in $r$.

{\bf Compact orientable surfaces of genus$>1$}:
A similar argument works for any genus $g>1$, with $P$ a polygon with $2g$ sides. Moreover, the fundamental group of a closed orientable surface of genus $g$ has uniformly exponential growth. Let $S$ be a system of generators for $\pi_1(\Sigma)$. Then it can be seen that $S$ contains some subset $A$ of $2g$ elements, which spans a subgroup of finite index in the abelianization ${\bf Z}^{2g}$ of $\pi_1(\Sigma)$. For an arbitrary $x\in A$, the group spanned by $A\setminus \{ x \}$ is of finite index in $\pi_1(\Sigma)$. Notice that it is the fundamental group of a noncompact surface. Hence it is a free group of rank $2g-1$, which has uniformly exponential growth. Therefore $\pi_1(\Sigma)$ also has uniformly exponential growth. A complete proof is available \cite{DeLaHarpe}. 

The result now follows in all cases from Lemma \ref{lemma:pi1growth}. \end{proof}

\subsection{Proof of Theorem \ref{Thm:3mfds-intractable-filters}}
\begin{proof}
 Once again, we will rely on Lemma \ref{lemma:pi1growth}, so our task first is to show that the manifolds in question have fundamental groups with exponential growth. 
 
 Let $Y$ be a compact connected orientable $3$--manifold, then the fundamental group of $Y$ has growth of exponential type, unless $Y$ is either homeomorphic to a geometric manifold modelled on one of the Thurston geometries ${\bf S}^3, {\bf S}^2\times {\bf R}, {\bf E}^3, {\bf Nil}$, or homeomorphic to a connected sum $L(2,1)\# L(2,1)$ of a lens space $L(2,1)$ whose fundamental group has order $2$, with itself. This follows, as in the case of surfaces, by the growth type of the possible fundamental groups. Assume $Y$ is not one of the manifolds just described. Then, by Thurston's geometrization programme \cite{Thurston}, completed by Perelman \cite{Perelman}, $\pi_1(Y)$ splits as a graph of groups. 
 
 If the splitting is not trivial, then at least one vertex group $G$ is isomorphic to a lattice in the isometry group of one of the geometries ${\bf H}^3, {\bf H}^{2} \times {\bf R}, {\bf Sol}$ or ${\rm PSL}_{2}({\bf R})$. Such a $G$ has a growth of exponential type, forcing $\pi_1(Y)$ to grow at least at the same rate and therefore also have exponential growth.
 
 If the splitting is trivial,$Y$ is a geometric manifold modelled on one of these geometries of exponential growth. Therefore, by Lemma \ref{lemma:pi1growth} in every one of these possible cases, the strategy to perform kernel translation along geodesics is not efficient. The only borderline case appears for the manifold $L(2,1)\# L(2,1)$, whose fundamental group does not grow exponentially \cite{BD}. Moreover, we can explicitly compute the order of the polynomial growth of the possible geometric manifolds. For ${\bf S}^3$ the growth is bounded, for ${\bf S}^2\times {\bf R}$ it is linear (because of the ${\bf R}$ factor), for ${\bf E}^3$ it is cubic, and for ${\bf Nil}$ it is at most quartic. 
  
 To verify this, we recall the definitions of these geometries. 
 
 The 3--sphere, with isometry group ${\rm SO}(4)$, here seen as the unit sphere in ${\bf R}^4$ with the induced metric. A family of manifolds modelled on ${\bf S}^3$ are the Lens Spaces $L(p,q)$. Here $p$ and $q$ are co-prime integers which define a $\mathbb{Z}_p$ action on $S^3 \subset \mathbb{C}^2$ by $(u, v)\mapsto (\omega^q u , \omega v )$ , where $\omega=e^{2i\pi p}$ and $L(p,q)$ is the quotient of ${\bf S}^3$ by this action. The growth of geodesics is thus bounded in this geometry.
 
 For ${\bf S}^2\times {\bf E}$, the isometry group of this geometry consists of the product of the spherical isometries times the isometries of the Euclidean line. There are only two orientable non-homeomorphic examples in this geometry. The product $S^{2}\times S^{1}$ and ${\bf R}{\rm P}^3 \#
 {\bf R}{\rm P}^3$, which is the only geometric 3--manifold that is also a non-trivial
 connected sum. On the $2$--sphere factor, the geodesics are bounded, and on the ${\bf E}$ factor, they contribute linear growth with respect to length. 
 
 For ${\bf E}^3$, Euclidean space, with isometry group ${\rm O}(3)\ltimes {\bf R}^3$, and the standard metric of ${\bf R}^3$, $ds^2=dx^2 + dy^2 + dz^2$. The 3-torus $T^3$ is modelled on this geometry. There are only six orientable and four non-orientable flat 3--manifolds, and their fundamental group classifies them. One of Bieberbach's theorems tells us that a group $\Gamma$ is isomorphic to a discrete group of isometries of ${\bf E}^n$ if and only if $\Gamma$ contains a subgroup of finite index that is free abelian of finite rank. Diffeomorphism classes of closed Euclidean $3$--manifolds are in one-to-one correspondence, via their fundamental group, with torsion-free groups containing a subgroup of finite index isomorphic to ${\bf Z}^3$. Therefore these flat manifolds are finitely covered by $T^3$, and hence their growth is at most polynomial.
 
 The Heisenberg group ${\bf Nil}^{3}$ can be defined as the following matrix group:
$${\bf Nil}^{3}=\left\{ \left(\begin{array}{ccc}
1 & x & y \\
0 & 1 & z \\
0 & 0 & 1 \end{array} \right) : x,y,z\in {\bf R} \right\}.$$
The isometry group $\is({\bf Nil}^{3})$ is the product of ${\bf Nil}^{3}$ by $S^{1}$ acting
as a group of automorphisms which rotate the $xy$-plane. We have an
exact sequence, 
$$ 0\rightarrow {\bf R}\rightarrow \is({\bf Nil}^{3})\rightarrow
\is ({\bf E}^{2})\rightarrow 0.$$

 Taking the subgroup $\Gamma$ of matrices with integer entries ${\bf Nil}^{3}$/$\Gamma$
is a circle fibration $\eta$ over the torus with orientable total space and 
 Euler number $e(\eta)=1$. In fact, any oriented circle bundle over a
 2--torus, which is not the 3--torus, has this kind of geometric structure.

To recognise the metric of ${\bf Nil}^{3}$ we identify it with ${\bf R}^{3}$,
$$ \left(\begin{array}{ccc}
1 & x & y \\
0 & 1 & z \\
0 & 0 & 1 \end{array} \right) \mapsto (x,y,z) .$$

Then ${\bf Nil}^{3}$ can be described as ${\bf R}^{3}$ with the metric $$ds^2=dx^2+dy^2+(dz-xdy)^2.$$ The growth rate of geodesics is proportional to a polynomial of order four. \end{proof}

\subsection{Proof of Theorem \ref{Thm:rat-hyp-intractable-filters}}
\begin{proof}
From Theorem \ref{thm:rat-hyp-h} we know that as $M$ is rationally hyperbolic, any smooth Riemannian metric on $M$ has $h>0$. So we may continue as in the proof of Theorem \ref{Thm:sface-intractable-filters}, we obtain that $C(x,y,\ell)$ grows faster than any polynomial in $\ell$ for any smooth Riemannian metric $g$ on $M$. 
Therefore, a procedure that requires averaging over the set of geodesics of length at most $\ell$ between points $x$ and $y$ in $M$ can not be performed in polynomial time. The result follows. 
\end{proof}

\subsection{Proof of Theorem \ref{Thm:htop-intractable-filters}}
\begin{proof}
We will now see how obstructions to an efficient computation of $C(x,y,\ell)$ emerge from the dynamics of the geodesic flow of $(M,g)$. Consider Ma\~n\'e's formula \ref{formula} for the topological entropy $h$ of the geodesic flow. 

Observe that, by the same reasoning as in the proof of Theorem \ref{Thm:sface-intractable-filters}, if an efficient algorithm for the computation of $C(x,y,\ell)$ exists, then $h$ must be zero. The result then follows from the deep work of Contreras---Theorems \ref{thm::Contreras-ThmC} \& \ref{thm::Contreras-ThmA} explained in section \ref{subsec:generic}--- whose main result implies the existence of an open and dense set of Riemannian metrics with positive topological entropy. \end{proof}

\section{Conclusions}\label{sec-conclusions}

\subsection{Related work}
The approach described here is related to particular G-CNNs \cite{KT18, CW16}, and these architectures have been implemented in various instances.
Some 3D voxel architectures are implicitly using the ideas presented here by considering their data as embedded in ${\bf R}^3$, these include VoxNet \cite{MS15}, V-Net \cite{MNA16}, Pointnet \cite{GGGOCA16}, and VoxelNet \cite{ZT18}. An architecture that incorporates normal directions to a voxel object and, in doing so, obtains better performance is NormalNet \cite{WCSBL19}. To reduce the memory footprint of voxel objects and improve performance, OctoNet incorporates sparsification of 3D data \cite{ROG17}. Moreover, submanifolds appear in a CNN architecture that has a linear cost for the number of active sites, with many computational economies while keeping state-of-the-art performance \cite{GEV18}.

Other proposed definitions of local convolution include defining polar geodesic coordinates patches \cite{MBBV15}, using spirals around points in meshes \cite{BBBZ19}, and using actions or a homogeneous space structure to define "correlations", instead of convolutions \cite{KT18, CBV18}. These strategies will be subject to the same topological constraints explained above, to the extent that they rely on averaging over geodesic paths between points.

As this note was being concluded, a helpful review and proposal for bringing more coherence to the growing body of work on geometric deep learning appeared \cite{BBCT}, as well as a comprehensive treatment of $G$--equivariant CNNs \cite{WFVW}.


\subsection{Perspectives}

We introduced a general definition of convolution for smooth manifolds. Furthermore, our methods provide theoretical foundations for some implementations already in use mentioned above. They can be defined in arbitrarily high dimensions and work globally for the manifold structure. In this way, we contribute to the development of the theoretical foundations of deep learning. 

We compared our proposal with some local convolution techniques and observed how group-theoretical properties show that they become computationally intractable as the genus grows in the case of surfaces. In short, it is undesirable to have to average over a group, as the group may not behave as well as needed for the analysis to be carried out. We also described topological obstructions from rational hyperbolicity and geometric/dynamical obstructions that originate in the growth of geodesic arcs between points. Moreover, this condition is generic in the space of Riemannian metrics. Finally, we pointed out how the finite metric space approximations first studied by Turing imply that an abelian structure must be used to define global and computable convolutions. 

The general framework presented here provides a rigorous and global definition of manifold convolution. It contributes to the theoretical foundations of higher dimensional convolutional architectures, beyond 3D.

\subsubsection*{Acknowledgments}

Thanks to Scott James, Armando Casta\~neda, Maks Ovsjanikov, Guido Mont\'ufar, and Minh Ha Quang for comments on this and earlier versions of this work. This research was supported in part by  DGAPA-UNAM PAPIIT grant IN104819. This project started when I visited the Institute for Pure and Applied Mathematics (IPAM), UCLA  in 2019, itself supported by NSF DMS-1440415. Thanks IPAM!

\providecommand{\bysame}{\leavevmode\hbox to3em{\hrulefill}\thinspace}

\end{document}